\documentclass[11pt, a4paper,notitlepage]{article}

\usepackage{amssymb}
\usepackage{amsmath}
\usepackage{amsthm}

\usepackage{epsfig}

\newcommand{\eps}{\varepsilon}

\newcommand{\FF}{\mathcal{F}}

\newcommand{\bx}{\mathbf{x}}
\newcommand{\bu}{\mathbf{u}}
\newcommand{\bv}{\mathbf{v}}
\newcommand{\bw}{\mathbf{w}}
\newcommand{\bq}{\mathbf{q}}
\newcommand{\bs}{\mathbf{s}}
\newcommand{\bk}{\mathbf{k}}
\newcommand{\bX}{\mathbf{X}}

\newcommand{\bz}{\mathbf{z}}

\newcommand{\half}{\frac{1}{2}}

\newtheorem{thm}{Theorem}[section]
\newtheorem{lem}[thm]{Lemma}

\newtheorem{cor}[thm]{Corollary}
\theoremstyle{definition}
\newtheorem{rem}[thm]{Remark}
\newtheorem{ex}[thm]{Example}

\newenvironment{prf}{\begin{proof}[\bf Proof]}{\end{proof}}

\def\infint{\int_{-\infty}^\infty}
\def\ex{{\rm E\,}}
\def\var{{\rm Var\,}}
\def\re{{\rm Re\,}}
\def\im{{\rm Im\,}}
\def\arg{{\rm Arg\,}}

\newcommand{\cov}{{\rm Cov}\,}
\newcommand{\E}{{\rm E}\,}
\newcommand{\txdi}{\tilde{X}^\Delta_i}
\newcommand{\xdi}{X^\Delta_i}

\begin{document}

\title{Nonparametric volatility density estimation
for discrete time models}

\author{Bert van Es, Peter Spreij \\
{\normalsize Korteweg-de Vries Institute for Mathematics}\\
{\normalsize Universiteit van Amsterdam}\\
{\normalsize Plantage Muidergracht 24}\\
{\normalsize 1018 TV Amsterdam}\\
{\normalsize The Netherlands} \and
Harry van Zanten\\
{\normalsize Division of Mathematics and Computer Science}\\
{\normalsize Vrije Universiteit}\\
{\normalsize De Boelelaan 1081 a}\\
{\normalsize 1081 HV Amsterdam}\\
{\normalsize The Netherlands}}

\maketitle

\begin{abstract}
We consider discrete time models for asset prices with a stationary volatility process. We aim at 
estimating the multivariate density of this process at a set of consecutive time 
instants.

A Fourier type deconvolution kernel density estimator based on the
logarithm of the squared process is proposed to estimate the
volatility density. Expansions of the bias and bounds on the
variance are derived.
\end{abstract}

\newpage

\section{Introduction}

Suppose that we have price data $S_0, S_1, \ldots$ of a certain
asset in a financial market. Let $X$ be the log-return process, defined by
$X_t=\log S_t-\log S_{t-1}$.
It is commonly believed that stochastic volatility models of the form
\begin{equation}\label{eq:s}
X_t = \sigma_t Z_t
\end{equation}
describe much of the observed behaviour of this type of data. Here
$Z$ is typically an i.i.d.\ noise sequence (often Gaussian) and
at each time $t$ the random variables $\sigma_t$ and $Z_t$ are independent.
We will assume that the process $\sigma$ is strictly stationary
and  that the (multivariate)
marginal distributions of $\sigma$  have a density
with respect to the Lebesgue measure on $(0,\infty)$.
Our aim is to construct a nonparametric estimator for the multivariate  density of 
$(\sigma_t,\ldots,\sigma_{t+p})$, and to study its
asymptotic behaviour.

Models that are used in the literature to
describe the volatility display rather different invariant
distributions. This observation lies at the basis of our point of view, which we pursue in this paper, that
nonparametric estimation procedures are by all means sensible tools to get
some insight in the behaviour of the volatility. Quite often in models that are used
in practice,
the invariant distributions of $\sigma$ are
unimodal. Since it is known that volatility clustering is an
often occurring phenomenon, it is hard to believe that this can
be explained by any of these models. Instead, one would expect in
such a case for instance  the distribution of $(\sigma_t,\sigma_{t+1})$ to have a density  that has concentration
regions around the diagonal with possibly peaks at
certain clusters of low and high volatility, a phenomenon that may lead to for instance bimodal
one-dimensional marginal distributions. Nonparametric density
estimation could perhaps reveal such a shape of the invariant
density of the volatility.

We will distinguish two classes of models in this paper. In both of them
we will assume that the noise sequence is standard Gaussian and  that  $\sigma$ is a strictly stationary,  positive
process satisfying a certain mixing condition. The way in which the bivariate 
process $(\sigma,Z)$, in particular its dependence structure, is further modelled differs however.
In the first class of models that we consider, we
assume that the process $\sigma$ is
predictable with respect to the filtration $\FF_t$ generated by the
process $Z$. 
Note that $\sigma_t$ is independent of $Z_t$ for each fixed time $t$.
We furthermore have that (assuming that the unconditional variances are finite)
$\sigma^2_t$ is equal to the conditional variance of $X_t$ given $\mathcal{F}_{t-1}$.
This class of models has become quite popular in
the econometrics literature. Financial data such as log-returns
of stock prices or exchange rates are believed to share a number
of stylized features, including for instance heavy-tailedness and
long-range dependence. Models of the type (\ref{eq:s}) have been
proposed to capture those features. A well-known family included
in the class (\ref{eq:s}) is the family of GARCH-models,
introduced by Bollerslev (1986). For the GARCH($p,q$)-model the
sequence $\{\sigma_t\}$ in (\ref{eq:s}) is assumed to satisfy
the equation
\begin{equation}\label{eq: garch}
\sigma^2_t = \alpha_0 + \sum_{i=1}^p \alpha_i X^2_{t-i} +
\sum_{j=1}^q \beta_j \sigma^2_{t-j},
\end{equation}
where the $\alpha_i$ and $\beta_j$ are nonnegative constants.
Under suitable assumptions, see Bougerol and Picard (1992), GARCH
processes are stationary and the statistical problem in this
case would be to estimate the coeficients $\alpha_i$ and
$\beta_j$ in~(\ref{eq: garch}).

In the second  class of models that we consider, we assume that the
whole process $\sigma$ is independent of the noise process $Z$. In this case, the natural underlying filtration
$\FF=\{\mathcal{F}_t\}_{t\geq 0}$ is
generated by the two processes $Z$ and $\sigma$ in the following way. For each $t$ the $\sigma$-algebra $\mathcal{F}_t$ is
generated by $Z_s$, $s\leq t$ and $\sigma_s$, $s\leq t+1$. This choice of the filtration
enforces $\sigma$ to be predictable. As in the first model the process $X$ becomes a martingale
difference sequence and we have again (assuming that the unconditional variances are finite)
that $\sigma^2_t$ is the conditional variance of $X_t$ given $\mathcal{F}_{t-1}$.
 An example of such a model is given in De Vries (1991), where
$\sigma$ is generated as an AR(1) process with $\alpha$-stable noise
($\alpha\in(0,1)$).

As we said before, we do not want to make a parametric assumption
such as (\ref{eq: garch}), but we still want to measure the
volatility of the data somehow. In the present paper we propose a
nonparametric statistical procedure for this problem.  Using ideas from
deconvolution theory, we will propose a procedure for the
estimation of the marginal density at a fixed point.  To assess the
quality of our procedure, we will 
derive expansions of the bias and bounds on the variance.
This will be done separately for the two kinds
of model classes outlined above.

\section{Primer on kernel type deconvolution}\label{primer}
\setcounter{equation}{0}

We briefly review the construction of the deconvolution kernel
density estimator  based on  i.i.d.\ observations, see also Wand and Jones (1995).
For simplicity we consider in this section the univariate case only. Recall that the
{\em characteristic function} or {\em Fourier transform} of a
density function $g$ is defined by
\begin{equation}\label{ft}
\phi_g(t)=\ex e^{itX} = \infint e^{itx}g(x)dx,
\end{equation}
where $X$ is a random variable with density function $g$. In the
standard deconvolution setting the random variable $X$ is equal to
the sum of two independent random variables, say $Y$, with unknown
density $f$, and $Z$, with known density $k$. So $g$ is the
convolution of $f$ and $k$ and
\begin{equation}\label{ft:2}
\phi_g(t)=\ex e^{itX} = \ex e^{it(Y+Z)}
 =  \ex e^{itY} \ex e^{itZ}
= \phi_f(t)\phi_k(t).
\end{equation}

The objective is to estimate $f$ from i.i.d. observations of
$X_1,\dots,X_n$ having density $g$. In identity (\ref{ft:2}) we
know $\phi_k(t)$ and we can estimate $\phi_g(t)$ by the
characteristic function of a {\em kernel estimator} $g_{nh}$ of
$g$. So
\begin{equation}
g_{nh}(x)={1\over n}\sum_{j=1}^n {1\over h}\,w\Big({{x-X_j}\over
h}\Big),
\end{equation}
where $w$ is an integrable function with integral one, called the
{\em kernel function}, and $h>0$ is a positive number, called the
{\em bandwidth}, governing the curvature of the estimate. The
kernel estimator itself is also a convolution of the empirical
distribution function $G_n$ of the observations and the rescaled
kernel function $w_h(x)=w(x/h)/h$. So, with $\phi_w$ the Fourier transform of $w$,
\begin{equation}
\phi_{g_{nh}}(t) = \phi_{w_h}(t) \infint e^{itx}dG_n(x)=
\phi_w(ht)\phi_{emp}(t),
\end{equation}
where
\begin{equation}
\phi_{emp}(t)=\infint e^{itx}dG_n(x)={1\over n}\sum_{j=1}^n
e^{itX_j}
\end{equation}
is called the {\em empirical characteristic function}. From
(\ref{ft:2}) we see that
\begin{equation}\label{ft:3}
{{\phi_w(ht)\phi_{emp}(t)}\over\phi_k(t)}
\end{equation}
is an obvious candidate to estimate $\phi_f$. Applying an inverse
Fourier transform we obtain an estimator of $f$. Define the
estimator $f_{nh}$ of $f$ as
\begin{equation}\label{fourest}
f_{nh}(x) = {1\over 2\pi}\infint
e^{-itx}{{\phi_w(ht)\phi_{emp}(t)}\over\phi_k(t)}\ dt.
\end{equation}

The inversion is allowed if the  function (\ref{ft:3}) is
integrable. In general this is not guaranteed. However, to
enforce integrability, we assume that $\phi_w$ has a bounded
support. Note that (\ref{fourest}) can be rewritten as
\begin{eqnarray}
\lefteqn{f_{nh}(x) ={1\over n}\sum_{j=1}^n {1\over 2\pi}\infint
{{\phi_w(ht)}\over\phi_k(t)}\ e^{-it(x-X_j)}dt}\nonumber\\
&=&{1\over nh}\sum_{j=1}^n {1\over 2\pi}\infint
{{\phi_w(s)}\over\phi_k(s/h)}\
e^{-is{(x-X_j)/h}}ds\label{sumfourest}\\ &=&{1\over
nh}\sum_{j=1}^n v_h\Big({{x-X_j}\over h}\Big),\nonumber
\end{eqnarray}
where
\begin{equation}\label{fourkernel}
v_h(x)={1\over 2\pi}\infint {{\phi_w(s)}\over\phi_k(s/h)}\
e^{-isx}ds.
\end{equation}
It is easy to see that the function $v_h$, and
hence the estimator $f_{nh}(x)$, is real valued.  Indeed, taking complex
conjugates, we get
\begin{eqnarray*}
\overline{v_h(x)} & = & {1\over 2\pi}\infint e^{isx} {\overline{{\phi_w(s)}}/\overline{\phi_k(s/h)}}\, ds
\\
& = & {1\over 2\pi}\infint  e^{-isx}{{\overline{\phi_w(-s)}}/\overline{\phi_k(-s/h)}}\,ds \\
& = & {1\over 2\pi}\infint  e^{-isx} {{\phi_w(s)}/\phi_k(s/h)}\,ds\\
& = & v_h(x).
\end{eqnarray*}

A popular performance measure for deconvolution kernel estimators
is the {\em mean squared error} (MSE). The MSE of $f_{nh}(x)$ is
defined as \mbox{$\ex(f_{nh}(x)-f(x))^2$}. To obtain asymptotic
expansions for the MSE, we need expansions for the bias and
variance of the estimator. The expectation of $f_{nh}(x)$ is
equal to the expectation of an ordinary kernel density estimator
of $f$ based on observations from $f$. We have
\begin{eqnarray*}
\lefteqn{\ex f_{nh}(x)=\int {1\over h}\,w\Big({{x-u}\over
h}\Big)f(u)du}\\
&=&f(x)+\tfrac{1}{2}h^2\int u^2w(u)du\,f''(x)+o(h^2),
\end{eqnarray*}
as $n\to \infty$, $h\to 0$  and $nh\to\infty$, provided that $w$
is symmetric
 and $f$ satisfies some smoothness
conditions, essentially twice differentiability at $x$. The
asymptotic variance of $f_{nh}(x)$ depends  on the tails of the
characteristic function of the density $k$. The smoother $k$, the
faster the tails of the characteristic function vanish and the
larger the asymptotic variance, see for instance Fan (1991).

\section{Construction of the estimators}\label{dtmodels}
\setcounter{equation}{0}

We consider the  model~(\ref{eq:s}), so
$X_t =\sigma_t Z_t.
$
If we
square this equation and take logarithms we get
\begin{equation}\label{eq: model2}
\log X^2_t  = \log\sigma^2_{t}+\log Z^2_t.
\end{equation}
Recall that under our assumptions for each $t$ the random variables $\sigma_t$ and
$Z_t$ are independent.
The density of $\log Z^2_t$, denoted by $k$, is given by
\begin{equation}\label{densityk}
k(x) = \frac{1}{\sqrt{2\pi}}\, e^{\tfrac{1}{2}x} e^{-
\tfrac{1}{2}e^{x}}.
\end{equation}
Its graph is given in Figure 1 below.

\begin{figure}[h]
$$ \epsfxsize=6cm\epsfysize=4cm\epsfbox{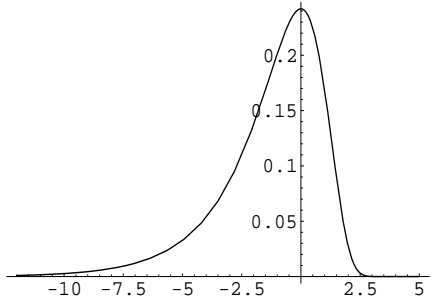} $$
\caption[]{The density function $k$ of $\log Z^2_t$.
\label{fig:1}}
\end{figure}

As in Section~\ref{primer}, it seems reasonable to use a
{\em deconvolution kernel density estimator} to estimate the
unknown density $f$ of $\log\sigma^2_{t}$. An estimate of the
density of $\sigma_t^2$ or $\sigma_t$ can then be obtained by a
simple transformation. Computing the characteristic function
$\phi_k$ of $\log Z^2$ we get, with $k(x)$ as in (\ref{densityk}),
\begin{equation}\label{phik}
\phi_k(t) = \infint e^{itx}k(x)dx =\tfrac{1}{\sqrt{\pi}}\,
2^{it}\, \Gamma(\tfrac{1}{2}+it),
\end{equation}
where the gamma function $\Gamma$ is defined for all complex $z$ with positive real
part by
\[
\Gamma(z)=\int_0^\infty t^{z-1}e^{-t}\,dt.
\]
The graphs of Re$(\phi_k)$, Im$(\phi_k)$ and $|\phi_k|$ are given
in Figures 2 and 3.

\begin{figure}[h]
$$ \epsfxsize=6cm\epsfysize=4cm\epsfbox{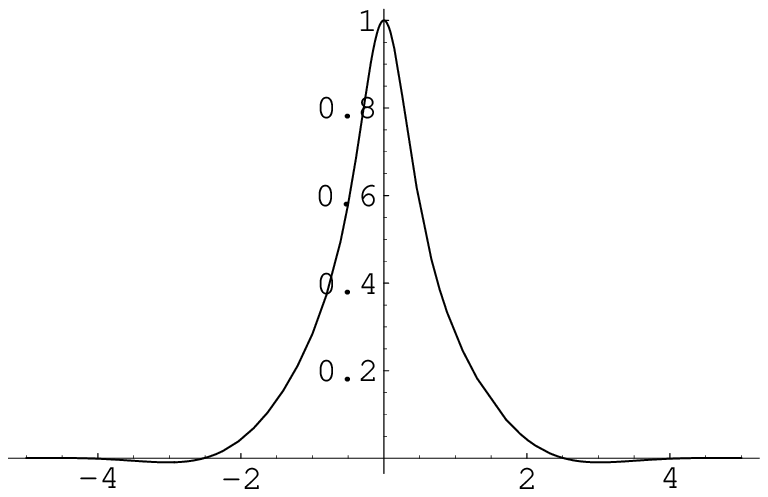}
\epsfxsize=6cm\epsfysize=4cm\epsfbox{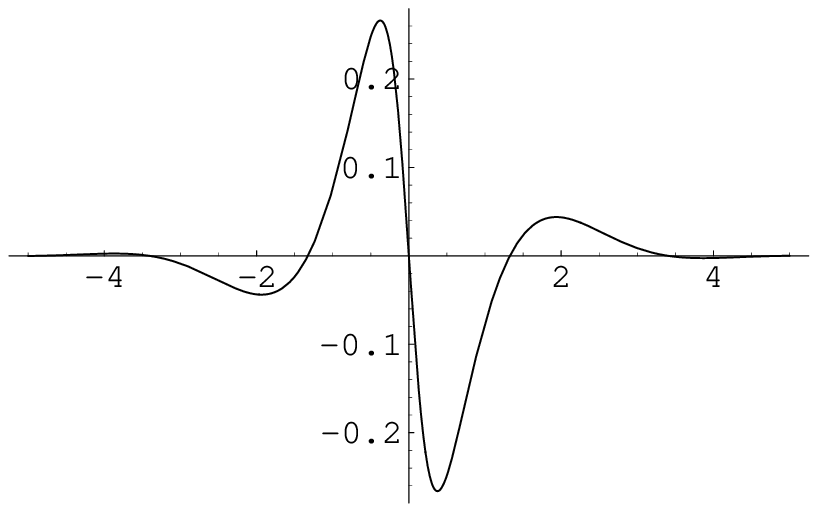} $$
\caption[]{The real and imaginary part of the characteristic
function $\phi_k$. \label{fig:2}}
\end{figure}
\begin{figure}[h]
$$ \epsfxsize=6cm\epsfysize=4cm\epsfbox{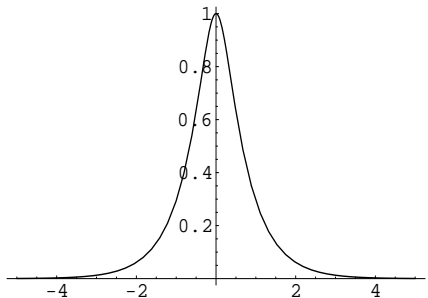} $$
\caption[]{The modulus of $\phi_k$. \label{fig:3}}
\end{figure}
\noindent

For the  model~(\ref{eq:s}) this leads to the
estimator
\begin{equation}\label{eq:fnh}
f_{nh}(x)={1\over nh}\sum_{j=1}^n v_h\Big({{x-\log X^2_j}\over
h}\Big)
\end{equation}
of the density $f$ of $\log \sigma^2_{t}$, with $v_h(x)$ as in
(\ref{fourkernel}). Note that, like in the previous section, this estimator is real valued.

The expression for the estimator of the density of the $p$-dimensional random vector $(\log \sigma^2_t,\ldots,\log
\sigma^2_{t-p+1})$ is similar.  We first introduce some auxiliary notation.
Let $p$ be fixed and write $\bx_j$ for a vector $(x_j,\ldots,x_{j-p+1})$. We use
similar boldface expressions for other (random) vectors. The kernel $\bw$ that we will use
in the multivariate case is just a product kernel, $\bw(\bx)=\prod_{j=1}^pw(x_j)$. 
Likewise $\bk(\bx)=\prod_{j=1}^p k(x_j)$.
Then with $\bv_h$ defined by
\begin{equation}\label{eq:bvh}
\bv_h(\bx)={1\over (2\pi)^p}\int_{\mathbb{R}^p} {{\phi_\bw(\bs)}\over\phi_\bk(\bs/h)}\
e^{-i\bs\cdot \bx}\,d\bs,
\end{equation}
where $\bs\in\mathbb{R}^p$ and $\cdot$ denotes inner product,
the multivariate density estimator is given by
\begin{equation}\label{eq:fnhp}
f_{nh}(\bx)={1\over (n-p+1)h^p}\sum_{j=p}^n \bv_h\Big({{\bx-\log \bX^2_j}\over
h}\Big),
\end{equation}
where we use $\log \bX^2_j$ to denote the vector $(\log X^2_j,\ldots,\log X^2_{j-p+1})$. 

\section{Asymptotics}
\setcounter{equation}{0}

The bias of the deconvolution estimator described in
Section~\ref{primer} will be seen to be the same as the bias of a kernel density
estimator based on independent observations from $f$. Hence, under standard smoothness assumptions, it is of order
$h^2$ as $h\to 0$. The variance of this type of deconvolution
estimator heavily depends on the rate of decay to zero of
$|\phi_k(t)|$ as $|t|\to\infty$. The faster the decay the larger
the asymptotic variance. In other words, the smoother $k$ the
harder the estimation problem. This follows for instance for
i.i.d. observations from results in Fan (1991) and for stationary
observations from the work of Masry (1991, 1993a,b).

The rate of decay of $|\phi_k(t)|$ for the density
(\ref{densityk}) is given by Lemma \ref{phiexpan} in Section
\ref{proofs}, where we show that
\begin{equation}
|\phi_k(t)|\sim\sqrt{2}\,e^{-\frac{1}{2}\pi |t|}, \quad \mbox{as}\ |t|\to \infty.
\end{equation}
By the similarity of the tail of this characteristic function to
the tail of a Cauchy characteristic function  we can expect the
same order of the mean squared error as in Cauchy deconvolution
problems, where it decreases logarithmically in $n$, cf. Fan
(1991) for results on i.i.d.\ observations. Note that this rate,
however slow, is faster than the one for normal deconvolution.

In the model~(\ref{eq: model2}) the sequence $\{\log X_t^2\}$ is not
independent, so results on the asymptotic behavior of the kernel
estimator of Section~\ref{primer} are not directly applicable. In
the literature also more general deconvolution problems have been
studied, where the  i.i.d.\ assumption has been relaxed. For
instance, the deconvolution model $X_j=Y_j+Z_j$, where
$\{Y_j,Z_j\}$ is a stationary sequence and the sequences $\{Z_j\}$
 $\{Y\}$ are independent has been treated by E. Masry (1991,
1993a,b).

Expansions for the variance of  the
deconvolution kernel estimator have been derived under several
mixing conditions.
Under the assumption that the volatility process is independent of the noise sequence,
the model (\ref{eq: model2})  fits into this scheme.
We will obtain similar results for the estimator when $\sigma$ (as a process) is not independent of $Z$,
but only predictable with respect to the filtration generated by $Z$.

Let us define the mixing conditions. For a certain process
$\{X_j\}$ let ${\cal F}_a^b$ be the $\sigma$-algebra of events
generated by the random variables $X_j,\ j=a,\dots,b$. Let the
mixing coefficient $\alpha_k$ be defined by
\begin{equation}\label{eq:alphamix}
\alpha_k=\sup_{A\in {\cal F}_{-\infty}^0,\ B\in {\cal
F}_k^\infty} |P(AB) -P(A)P(B)|.
\end{equation}
We call a process $\{X_j\}$ {\em  strongly mixing} if
$\alpha_k\to 0$ as $k\to\infty$.

To obtain expansions for the bias and variance we also need
conditions on the kernel function $w$ such as bounded support of
its characteristic function $\phi_w(t)$. Moreover, the rate of
decay to zero of $\phi_w(t)$ at the boundary of its support turns
up in the asymptotics. The complete list of assumptions on $w$
that we use is the following.

\bigskip

\noindent{\bf Condition W.} Let $w$ be a real symmetric function
satisfying
\begin{enumerate}
\item
$\infint |w(u)|du < \infty,$
\item
$\infint w(u)du=1,$
\item
$\infint u^2w(u)du<\infty,$
\item
$\lim_{|u|\to\infty} w(u)=0$,
\item
$\phi_w$, the characteristic function of  $w$ has support [-1,1],
\item
$\phi_w(1-t)=At^\alpha + o(t^\alpha),\quad\mbox{as}\ t\downarrow
0 $ for some $\alpha > 0$.
\end{enumerate}

\bigskip
\noindent Note that by Fourier inversion  these conditions imply
that $w$ is bounded and Lipschitz. More precisely, we have
\begin{equation}\label{wbounds}
|w(x)|\leq {1\over 2\pi}\quad\mbox{and}\quad |w(x+u)-w(x)|\leq
{1\over 2\pi}|u|.
\end{equation}
An example of such a kernel, from Wand (1998), with $\alpha=3$
and $A=8$, is
\begin{equation}
w(x)={{48x(x^2-15)\cos x - 144(2x^2-5)\sin x}\over{\pi x^7}}.
\end{equation}
It has characteristic function
\begin{equation}
\phi_w(t)=(1-t^2)^3,\quad |t|\leq 1.
\end{equation}
The next theorem, whose proof can be found in Section~\ref{proofs}, establishes the expansion of the bias and an order
bound on the variance of our estimator under a strong mixing
condition. Under broad conditions this mixing condition is
satisfied if the process $\sigma$ is a Markov chain, since then
convergence of $\alpha_k$ to zero takes place at an  {\em
exponential rate}, see Theorems 4.2 and Theorem 4.3 of Bradley
(1985)
 for precise statements. Similar behaviour occurs for ARMA processes
with absolutely
 continuous distributions of the noise terms (Bradley (1985), Example
6.1).

\begin{thm}\label{discrasthm}
Assume that the process $X$ is strongly mixing with
coefficient $\alpha_k$ satisfying
$$
\sum_{j=1}^\infty \alpha_j^{\beta}<\infty,
$$
for some $\beta\in (0,1)$. Let the kernel function $w$ satisfy
Condition W and let the density $f$ of the $p$-vector $(\log\sigma^2_1,\ldots,\log \sigma^2_p)$ be
bounded and twice continuously differentiable with bounded second order partial derivatives.
Assume that $\sigma$ is a predictable process
with respect to the filtration generated by the process $Z$. 
Then we have for the estimator of the multivariate density defined as
in~(\ref{eq:fnhp}) and $h\to 0$
\begin{equation}\label{discrasthm:1p}
\ex f_{nh}(\bx)= f(\bx)+\tfrac{1}{2}h^2\int \bu^\top \nabla^2f(\bx)\bu\, \bw(\bu)\,d\bu+o(h^2)
\end{equation}
and
\begin{equation}\label{discrasthm:2p}
\var f_{nh}(\bx) = O\big(\tfrac{1}{n}\,(h^{2\alpha- \beta}\,e^{\pi/h})^p\big).
\end{equation}
\end{thm}

\bigskip

\begin{thm}\label{discrasthmp}
Assume that the process $\sigma$ is strongly mixing with
coefficient $\alpha_k$ satisfying
$$
\sum_{j=1}^\infty \alpha_j^{\beta}<\infty,
$$
for some $\beta\in (0,1)$. Let the kernel function $w$ satisfy
Condition W and let the density $f$ of the $p$-vector $(\log\sigma^2_1,\ldots,\log \sigma^2_p)$ be
bounded and twice continuously differentiable with bounded second order partial derivatives.
Assume furthermore that $\sigma$ and $Z$
are  independent processes. Then
 the multivariate density estimator $f_{nh}$ satisfies the same bias expansion as
in Theorem~\ref{discrasthm}. For the variance we have the sharper bound
\begin{equation}\label{discrasthm:3p}
\var f_{nh}(\bx) = O\big(\tfrac{1}{n}\,(h^{2\alpha}\,e^{\pi/h})^p\big).
\end{equation}
\end{thm}

\bigskip

\begin{rem}
Because of the exponential factor in the variance bound, in order
to obtain consistency, one has to take essentially $h\geq \pi/\log
n$, see also Stefanski (1990) for a related problem. On the other
hand we would like to minimize the bias, so the choice $h=
\pi/\log n$ is optimal. Both bias and variance decay at a
logarithmic rate for this choice of bandwidth. This seems
disappointing, however Fan (1991) shows for the i.i.d. situation
of Section 2 that we can not expect anything better.
\end{rem}

\begin{rem}
Notice that the results in Masry (1993a,b) establishing strong
consistency, rates of convergence and asymptotic normality are
not useful here, because the condition that $\phi_k$ has either
purely real or purely imaginary tails is not satisfied.
\end{rem}

\begin{rem}
Note that our assumptions in Theorem~\ref{discrasthm} are slightly different from those of Masry (1991). One of
the essential facts that are used in the proof is the mixing property of $X$. If
$\sigma$ and $Z$ are independent processes this is implied by a similar assumption
on the $\sigma$ process itself as in Masry (1991).
\end{rem}

\begin{rem}\label{rem:ftilde}
In the case where the processes $\sigma$ and $Z$ are independent,
the estimators $f_{nh}(x)$ have the following property.
\begin{equation}\label{eq:fff}
\tilde{f}_{nh}(x):=\ex [f_{nh}(x)|\mathcal{F}^\sigma]=\frac{1}{nh}\sum_{j=1}^n
w\big(\frac{x-\log\sigma^2_j}{h}\big),
\end{equation}
where $\mathcal{F}^\sigma$ denotes the $\sigma$-algebra generated by the whole
process $\sigma$.
Thus the $\tilde{f}_{nh}(x)$ would be ordinary kernel density estimators, if the
$\sigma^2_j$ could be observed.

Equation~(\ref{eq:fff}) is seen to be true as follows.
Write $u_j=\log X_j^2$ and use similar
notation for $\zeta_j=\log Z_j^2$ and $\tau_j=\log\sigma^2_j$.
Then
\begin{eqnarray*}
\ex [v_h(\frac{x-u_j}{h})|\mathcal{F}^\sigma] & = &
\frac{1}{2\pi}\int \ex e^{is\zeta_j/h}\frac{\phi_w(s)}{\phi_k(s/h)}
e^{-is(x-\tau_j)/h}\, ds \\
& = & \frac{1}{2\pi}\int \phi_w(s)
e^{-is(x-\tau_j)/h}\, ds \\
& = & w\big(\frac{x-\tau_i}{h}\big).
\end{eqnarray*}
The result now follows. Of course, the analogous statement for the multivariate
density estimator is equally true. One has
\begin{equation}\label{eq:fffp}
\tilde{f}_{nh}(\bx):=\ex [f_{nh}(\bx)|\mathcal{F}^\sigma]=\frac{1}{nh^p}\sum_{j=p}^n
w\big(\frac{\bx-(\log\sigma^2_j,\ldots,\log\sigma^2_{j-p+1})}{h}\big),
\end{equation}
\end{rem}

\begin{rem}
Better bounds on the asymptotic variance than in Theorem~\ref{discrasthm} can be obtained under
stronger mixing conditions. Consider for instance  {\em uniform
mixing}. In this case the mixing coefficient $\phi_t$ is defined
for $t >0$ as
\begin{equation}\label{eq:betamix}
\phi_t =  \sup_{A\in {\cal F}_{-\infty}^0, B\in {\cal F}_t^{\infty}}
|P(A|B) -P(A)|.
\end{equation}
Similar to strong mixing,  a process is called  {\em uniform
mixing} if $\phi_t\to 0$ for $t\to \infty$. Obviously,  uniform
mixing implies strong mixing. As a matter of fact, one has the
relation
\[
\alpha_t \leq \tfrac{1}{2}\phi_t.
\]
See Doukhan (1994) for this inequality and many other mixing
properties. If $\{\sigma_t\}$ is uniform mixing with coefficient
 $\phi$ satisfying
$\sum_{j=1}^\infty \phi(j)^{1/2}<\infty$, then  the variance
bound~(\ref{discrasthm:2p}) can be replaced with
\begin{equation}\label{contasthm:3p}
\var f_{nh}(\bx) = O\Big({1\over
n}\,(h^{2\alpha}e^{\pi/h})^p\Big).
\end{equation}
The proof of the latter bound runs similarly to the strong-mixing
bound as given in section~\ref{proofs}. The essential difference is that in
equation~(\ref{eq:mnh})
we use Theorem 17.2.3 of Ibragimov and Linnik (1971) with
$\tau=0$  instead of Deo's (1973) lemma, as in the proof of
Theorem 2 in Masry (1983). The result is that we can now bound the term $M_{nh}$ of equation (\ref{eq:mnh}) by
a constant times $\sum_{j=1}^{n-p+1} \varphi_j^{1/2} \ex W_0^2$. After this step the
proof is essentially unchanged. Use the estimate $\ex W_0^2 \leq C h^p ||v||_2^2$ to
finish the proof. Notice that this bound on the variance is of the same order as the
one we obtained in Theorem~\ref{discrasthmp}, where $\sigma$ was only assumed to be strongly mixing.
This bound cannot be improved upon by strengthening the assumption to uniform
mixing.
\end{rem}

\begin{rem}
An example of an observed process that is stongly mixing and that belong to the first
model class is a GARCH$(p,q)$ process. It has been shown in Carasso and Chen (2002) (see also Boussama (1998))
that such a process is $\beta$-mixing with exponentially decaying $\beta$-mixing coefficients. Hence this
process is also $\alpha$-mixing,
since the $\beta$-mixing coefficient $\beta_k=\ex\, \mbox{ess sup} \{
|P(A|\mathcal{F}_k^\infty)-P(A)|:A\in\mathcal{F}^0_{-\infty} \}$
satisfies the inequality $2\alpha_k\leq \beta_k$ (see Doukhan
(1994)). Notice that we also have that the assumption of Theorem~\ref{discrasthm} on
the $\alpha$'s is satisfied in this case.
\end{rem}

\section{Proofs}\label{proofs}
\setcounter{equation}{0}

All the estimators that we proposed involve the functions
$\phi_k$ and $\phi_w$. For these functions and related ones we
need expansions and order estimates. These are collected in the
lemmas of this subsection.
\begin{lem}\label{phiexpan}
 For $|t|\to\infty$ we have
\begin{eqnarray*}\label{phikexp}
|\phi_k(t)|\ \ &=&\sqrt{2}\,e^{-\frac{1}{2}\pi
|t|}(1+O(\tfrac{1}{|t|})),\\
\re
\phi_k(t)&=&|\phi_k(t)|[\cos(t\log(\sqrt{1+4t^2}-t))+O(\tfrac{1}
{|t|})],\\
\im
\phi_k(t)&=&|\phi_k(t)|[\sin(t\log(\sqrt{1+4t^2}-t))+O(\tfrac{1}
{|t|})].
\end{eqnarray*}
\end{lem}

\begin{prf}
By the Stirling formula for the complex gamma function, cf.
Abramowitz and Stegun (1964) Chapter 6, we have
\begin{equation}
\log \Gamma(z) = (z-\tfrac{1}{2})\log z -z + \tfrac{1}{2} \log
2\pi + O(\tfrac{1}{|z|}),
\end{equation}
as $|z|\to\infty$ and $|\arg z |< \pi $ for some $\delta >0$. So
for $z=\frac{1}{2}+it$ and $|t|\to\infty$ we get
\begin{eqnarray*}
\lefteqn{\log \Gamma(\tfrac{1}{2}+it) = it \log( \tfrac{1}{2}+it)
-(\tfrac{1}{2}+it) + \tfrac{1}{2} \log 2\pi +
O(\tfrac{1}{|t|})}\\ & =& it(\log|\tfrac{1}{2}+it| + i \arg(
\tfrac{1}{2}+it)) -(\tfrac{1}{2}+it) + \tfrac{1}{2} \log 2\pi +
O(\tfrac{1}
{|t|})\\
&=&-t\arg(\tfrac{1}{2}+it))-\tfrac{1}{2}+ \tfrac{1}{2} \log 2\pi
 +i(t\log
|\tfrac{1}{2}+it|-t) +O(\tfrac{1}{|t|}).
\end{eqnarray*}
Taking the modulus of the exponent the imaginary part  vanishes
and  we get
\begin{eqnarray*}
\lefteqn{ |\Gamma(\tfrac{1}{2}+it)|=
\exp(-t\arg(\tfrac{1}{2}+it))-\tfrac{1}{2}+ \tfrac{1}{2} \log
2\pi+O(\tfrac{1}{|t|}))}\\ &
=&\sqrt{2\pi}\exp(-t\arctan{2t}-\tfrac{1}{2}+O(\tfrac{1}{
|t|}))\\ &
=& \sqrt{2\pi}\exp(-\tfrac{1}{2}\pi|t| +O(\tfrac{1}{|t|}))\\
&=& \sqrt{2\pi}\exp(-\tfrac{1}{2}\pi|t|)(1 +O(\tfrac{1}{|t|})).
\end{eqnarray*}
Here we have used the expansion
$t\arctan{t}=t(\frac{1}{2}\pi-\arctan(1/t))=
\frac{1}{2}\pi t-1 +O(1/t)$, as
$t$ tends to infinity. For negative $t$ a similar expansion
holds. Since $2^{it}=\exp(it\log 2)$ has modulus one,
substituting this expansion in (\ref{phik}) now proves the first
statement of the lemma. The argument of $\Gamma(\frac{1}{2}+it)$
satisfies
\begin{eqnarray*}
\lefteqn{\arg(\Gamma(\tfrac{1}{2}+it))=
t\log|\tfrac{1}{2}+it|-t+O(\tfrac{1}{|t|})}\\ & =& -t\log 2 +
t\log|1+2it|-t+O(\tfrac{1}{|t|}) .
\end{eqnarray*}
So, since $\arg(2^{it})=t\log 2$, we have
\begin{equation*}
\arg( \phi_k(t))=t\log(\sqrt{1+4t^2}-t)+O(\tfrac{1}{ |t|}),
\end{equation*}
which proves the second and third statement of the lemma.
\end{prf}
Consider now the function $v_h$ defined in~(\ref{fourkernel}).
\begin{lem}\label{l2exp} We have the following order estimate for the
$L^2$ norm of $v_h$. For $h\to 0$
\begin{equation}
\|v_h\|_2= O(h^{\half+\alpha}e^{\pi/2h}).
\end{equation}
\end {lem}

\begin{prf}
By Parseval's identity
$$
\|v_h\|^2= {1\over 2\pi}\int_{-1}^1
\big|{{\phi_w(s)}\over\phi_k(s/h)}\Big|^2dx.
$$
Write
\begin{eqnarray}
\lefteqn{\int_{-1}^1\Big|{{\phi_w(s)}\over\phi_k(s/h)}\Big|^2ds
}\nonumber\\
&\leq&
 \tfrac{1}{2}\int_{-1}^1|\phi_w(s)|^2 \,e^{\pi
|s/h|}ds \label{een}\\
&+& \int_{-1}^1|\phi_w(s)|^2\Big|{1\over|\phi_k(s/h)|^2}-
\tfrac{1}{2}e^{ \pi |s/h|}\Big|ds \label{twee}
\end{eqnarray}
The integral in (\ref{een}) can be rewritten as
\begin{eqnarray*}
\lefteqn{\int_{-1}^1|\phi_w(s)|^2 \,e^{\pi |s/h|}ds}\\
&=& e^{ \pi/h}\int_{-1}^1|\phi_w(s)|^2 \,e^{ \pi
(|s/h|-(1/h))}ds\\
&=& 2e^{ \pi/h}\int_{0}^1|\phi_w(s)|^2 \,e^{ \pi
(|s/h|-(1/h))}ds\\
&=& 2e^{ \pi/h}h\int_{0}^{1/h}|\phi_w(1-hv)|^2
\,e^{ \pi ((1-hv)/h-(1/h))}dv\\
&=& 2e^{ \pi/h}h^{1+2\alpha}\int_{0}^{1/h}\Big|{\phi_w(1-hv)
\over  (hv)^{\alpha}}
\Big|^2v^{2\alpha} \,e^{- \pi v}dv\\
&\sim&
2e^{ \pi/h}h^{1+2\alpha}A^2\int_0^{\infty}v^{2\alpha} e^{- \pi v}dv\\
&=& 2e^{ \pi/h}h^{1+2\alpha}(
\pi)^{-1-2\alpha}A^2\Gamma(2\alpha+1),
\end{eqnarray*}
by the dominated convergence theorem. Omitting constants, we can
rewrite the integral  (\ref{twee}) as
\begin{eqnarray*}
\lefteqn{\int_{-1}^1|\phi_w(s)|^2e^{ \pi |s/h|} \Big|{ 2 e^{- \pi
|s/h|}\over|\phi_k(s/h)|^2}- 1\Big|ds}\\ &=& e^{
\pi/h}\int_{-1}^1|\phi_w(s)|^2 \Big|{ 2 e^{- \pi
|s/h|}\over|\phi_k(s/h)|^2}- 1\Big|
e^{ \pi (|s/h|-(1/h))}ds\\
&=& 2e^{ \pi/h}\int_0^1|\phi_w(s)|^2 \Big|{ 2 e^{- \pi
|s/h|}\over|\phi_k(s/h)|^2}- 1\Big|
e^{ \pi (|s/h|-(1/h))}ds\\
&=& 2h^{1+2\alpha}e^{\pi/h}\int_0^{1/h}\Big|{{|\phi_w(1-hv)|}
\over (hv)^\alpha}\Big|^2 \Big|{ 2 e^{- \pi
(1/h-v)}\over|\phi_k(1/h-v)|^2}- 1\Big|v^{2\alpha} e^{- \pi
v}dv\\
&=& 2h^{1+2\alpha}e^{\pi/h}o(1),
\end{eqnarray*}
by the dominated convergence theorem. We have used the fact that
both the functions $\phi_w(1-u)/u^\alpha$ and (see Lemma~\ref{phiexpan}) $|(2\exp(-\pi
u)/|\phi_k(u)|^2)- 1|$ are bounded and that the second function
is of order $O(1/u)$ as $u$ tends to infinity. This shows that
the   term   (\ref{twee}) is negligible with respect to
(\ref{een}).
\end{prf}

\begin{cor}\label{cor:vp}
 The $L^2$-norm of the function $\bv_h$, defined in (\ref{eq:bvh}) is of order
$O\big(h^{p(\half+\alpha)}e^{p\pi/2h}\big)$.
\end{cor}

\begin{prf}
This follows from the product form of $\bv_h$ given by
$\bv_h(\bs)=\prod_{j=1}^p v(s_j)$.
\end{prf}

\noindent {\bf Proof of Theorem~\ref{discrasthm}.} The
expansion~(\ref{discrasthm:1p}) follows from Theorem 1 in Masry (1991). To
prove the variance bound (\ref{discrasthm:2p}) we argue as in the proof
of Theorem 2 in the same paper. First we give a bound on the variance in terms of
the $L_2$-norm of the function $\bv_h$ and then we exploit the asymptotic expansion of
the characteristic function $\phi_k$ as given in Lemma~\ref{phiexpan} to get a
sharper bound on the $L_2$-norm of $\bv_h$ than Masry in his Proposition 3 by taking the behaviour of $\phi_w$ at the
boundary of its support into account. Some details follow.

Argueing as in Masry (1991) we can show that
\[
\var f_{nh}(\bx)=O\big(\frac{||\bv_h||^2_2}{nh^p} + M_{nh} \big),
\]
with (up to a multiplicative constant)
\begin{equation}\label{eq:mnh}
M_{nh}=\frac{1}{nh^{2p}}\sum_{j=p}^{n}\cov(W_j,W_0),
\end{equation}
where $W_j=\bv_h(\frac{\bx-\log \bX_j}{h})$.

Applying a lemma by Deo (1976), we can bound for strong mixing process $X$ with
mixing coefficients $\alpha_j$ the term $M_{nh}$ by a constant (not depending on $n$ and $h$) times
\[
\frac{1}{nh^{2p}}\sum_{j=p}^n \alpha_{j-p+1}^\beta \{\ex |W_j|^{2/(1-\beta)}\ex |W_0|^{2/(1-\beta)}\}^{(1-\beta)/2},
\]
which, by stationarity, becomes
\[
\frac{1}{nh^{2p}}\sum_{j=p}^n \alpha_{j-p+1}^\beta (\ex |W_0|^{2/(1-\beta)})^{1-\beta}.
\]
Observe now that, by boundedness of the density of $\log X_j^2$,
the term $\ex |W_0|^{2/(1-\beta)}$ can be bounded by a constant times
$h^p||\bv_h||_{2/(1-\beta)}^{2/(1-\beta)}$ and that we can therefore write
\[
\var f_{nh}(\bx)=O\big(\frac{||\bv_h||^2_2}{nh^p} +
\sum_{j=p}^{n}\alpha_j^{\beta}
\frac{||\bv_h||_{2/(1-\beta)}^2}{nh^{p(1+\beta)}}\big).
\]
The proof will be finished by application of Corollary~\ref{cor:vp}, which gives the
$L^2$-norm of $\bv_h$, and an estimate of the $L^{2/(1-\beta)}$-norm of $\bv_h$. For the latter
one we have the inequalities $||\bv_h||_{2/(1-\beta)}\leq ||\bv_h||_\infty^\beta ||\bv_h||_2^{1-\beta}$ and
$||\bv_h||_\infty \leq C ||\bv_h||_2$ for some constant $C$ by the fact that $\phi_w$
has compact support. As a result we get $||\bv_h||_{2/(1-\beta)}\leq C ||\bv_h||_2$ and
that $M_{nh}$ is less than a constant times $||\bv_h||_2/nh^{p(1+\beta)}$. The bound
on $\var f_{nh}(\bx)$ of theorem~\ref{discrasthm} now follows.
\hfill$\square$
\medskip\\
{\bf Proof of Theorem~\ref{discrasthmp}}.
Let $\mathcal{F}^\sigma$ be the $\sigma$-algebra generated by the
process $\sigma$. We use the decomposition
\begin{equation}\label{eq:decvar}
\var f_{nh}(\bx) = \ex\var(f_{nh}(\bx)|\mathcal{F}^\sigma) + \var\tilde{f}_{nh}(\bx),
\end{equation}
with $\tilde{f}_{nh}(\bx)$ as in Remark~\ref{rem:ftilde}. We now consider the first
term in~(\ref{eq:decvar}).
Let $\bz_j=(\log Z_j^2,\ldots,\log Z^2_{j-p+1})$ and
$\bq_j=(\log\sigma^2_j,\ldots,\log\sigma^2_{j-p+1})$.
Since the $Z_i$ are independent given
$\mathcal{F}^\sigma$ we can bound the conditional variance by
\[
\frac{1}{n^2h^{2p}}\sum_{j=p}^n \ex[ \big(\bv_h(\frac{\bx-\bq_j-\bz_j}{h})\big)^2|\mathcal{F}^\sigma]
\]
which is by conditional independence and stationarity equal to
\[
\frac{1}{nh^{2p}}\int \big(\bv_h(\frac{\bx-\bq_0-\bz}{h})\big)^2k(\bz)\,d\bz
\leq \frac{C}{nh^p}||\bv_h||^2_2,
\]
with $C$ the maximum of $\bk$, the density of $\bz_0$. Therefore the first term in
~(\ref{eq:decvar}) is of order $||\bv_h||_2^2/nh^p$, so of order $O\big(h^{p(1+2\alpha)}e^{p\pi/h}/nh^p\big)$. 

The second term of (\ref{eq:decvar}) is treated next. We have with
$U_j=\bw(\frac{\bx-\bq_j}{h})$
\begin{eqnarray*}
\var\tilde{f}_{nh}(\bx) & = & \frac{1}{n^2h^{2p}}\sum_{j} \var U_j
+\frac{2}{n^2h^{2p}}\sum_{i<j}\cov (U_i,U_j).
\end{eqnarray*}
The first term reduces by stationarity to
$\frac{1}{nh^{2p}}\var U_1$ which can be bounded by a constant times $||\bw||_2^2/nh^p$,
since $(\log\sigma_1^2,\ldots,\log\sigma_p^2)$ has by assumption a bounded density. For the second term we
proceed as in the proof of Theorem~\ref{discrasthm}. Using stationarity we write it as
\[
\frac{2}{n^2h^{2p}}\sum_{k=1}^n (n-k)\cov (U_k,U_0).
\]
We split the summation into two parts. In the first part we consider
\[
\sum_{k=1}^{p-1} (n-k)\cov (U_k,U_0).
\]
whose absolute value can be bounded in view of the Cauchy-Schwarz inequality and
stationarity by $(p-1)n\ex U_0^2$, which is bounded by $(p-1)nh^p||\bw||_2^2$.

The absolute value of the second part
\[
\sum_{k=p}^n (n-k)\cov (U_k,U_0)
\]
can be bounded by invoking once more Deo's result by
\[
n\sum_{k=p}^n\alpha_{k-p+1}^\beta\big(\ex|U_0|^{2/(1-\beta)}\big)^{1-\beta},
\]
which is less than
\[
nh^{p(1-\beta)}||w||_{2/(1-\beta)}^2\sum_k\alpha_{k-p+1}^\beta .
\]
Hence we have that $\var \tilde{f}_{nh}(\bx)$ is of order $1/nh^{p(1+\beta)}$.

Combining the obtained order estimates for the two terms of~(\ref{eq:decvar}) and
using the $L^2$-norm of the function $\bv_h$ gives the desired result.
\hfill$\square$

\section*{References}

\begin{verse}

Abramowitz, M.\ and Stegun, I.\ (1964), {\em Handbook of
Mathematical Functions, ninth edition}, Dover, New York.

Bollerslev, T.\  (1986), Generalized autoregressive conditional
heteroscedasticity,  {\em J.\ Econometrics} {\bf 31}, 307--321.

Bougerol, P.\  and N.\ Picard (1992a), Strict stationarity of
generalized autoregressive processes, {\em Ann.\ Probab.} {\bf
20}, no.~4, 1714--1730.

Bougerol, P.\  and N.\ Picard (1992b), Stationarity of GARCH
processes and of some nonnegative time series, {\em J.\
Econometrics} {\bf 52}, no.~1-2, 115--127.

Boussama, F.\ (1998), {\em Ergodicit\'e, m\'elange et estimation dans les mod\`eles
GARCH}, PhD, thesis, Universit\'e Paris 7.

Bradley, R.C.\ (1985), Basic properties of strong mixing
conditions, in {\em Dependence in Probability and Statistics, E.\
Eberlein and M.S.\ Taqqu Eds.}, Birkha\"user.

Carasso, M.\  and  Chen, X.\ (2002) Mixing and moment properties of various GARCH and
stochastic volatility models, {\em Econometric Theory} {\bf 18}, 17--39.

Deo, C.M.\ (1973), A note on empirical processes for strong
mixing processes, {\em Ann. Probab.} {\bf 1}, 870--875.

Doukhan, P.\ (1994), {\em Mixing, Properties and Examples},
Springer-Verlag.

Fan, J.\ (1991), On the optimal rates of convergence for
nonparametric deconvolution problems, {\em Ann. Statist.} {\bf
19}, 1257--1272.


Hewitt, E.\ and Stromberg K.\ (1965), {\em Real and Abstract
Analysis}, Springer Verlag, New York.

Ibragimov, I.A., and Linnik (1971), {\em Independent and
stationary sequences of random variables}, Wolters-Noordhoff.

Masry, E.\ (1983), Probability density estimation from sampled
data, {\em IEEE Trans. Inform. Theory}  {\bf 29}, 696--709.

Masry, E.\ (1991), Multivariate probability density deconvolution
for stationary stochastic processes, {\em IEEE Trans. Inform.
Theory}  {\bf 37}, 1105--1115.

Masry, E.\ (1993a), Asymptotic normality for deconvolution
estimators of multivariate densities of stationary processes,
{\em J. Mult. Anal.}  {\bf 44}, 47--68.

Masry, E.\ (1993b), Strong consistency and rates for
deconvolution of multivariate densities of stationary processes,
{\em Stoc. Proc. and Appl.} {\bf 475}, 53--74.

Stefanski, L.A.\ (1990), Rates of convergence of some estimators in a class of
deconvolution problems, {\em Statist.\ Probab.\ Lett.} {\bf 9}, 229--235.

De Vries, C.G.\ (1991), On the relation between GARCH and stable processes, {\em J.
Econometrics} {\bf 48}, 313--324.

Wand, M.P..\ (1998), Finite sample performance of deconvolving
kernel density estimators, {\em Statist.\ Probab.\ Lett.} {\bf 37},
131--139.

Wand, M.P.\ and Jones, M.C.\ (1995), Kernel Smoothing, Chapman and Hall,
London.

\end{verse}

\end{document}